\documentclass[11pt,leqno]{article}
\usepackage{a4}
\usepackage[T1]{fontenc}
\usepackage[english]{babel}
\usepackage{latexsym}
\usepackage{amsmath}
\usepackage{amssymb}
\usepackage{mathrsfs}
\usepackage{fig4tex184}
\usepackage{color}
\usepackage{cite}
\usepackage{titling}
\usepackage{bbold}
\makeatletter
\@addtoreset{equation}{section} 
\makeatother
\definecolor{move}{rgb}{.3,.1,.8}
\newcommand{\ud}{\mathrm{d}}
\newcommand{\vfi}{\varphi}

\newcommand{\foo}{\vfi_{1,1}}
\newcommand{\phot}{\vfi_{1,2}}
\newcommand{\phto}{\vfi_{2,1}}
\newcommand{\ftt}{\vfi_{2,2}}

\newcommand{\supp}{\mathrm{supp}}
\newcommand{\re}{\mathrm{Re}}
\newcommand{\im}{\mathrm{Im}}
\newcommand{\e}{\mathrm{e}}
\newcommand{\id}{\mathrm{Id}}
\newcommand{\Rn}{\mathbb{R}^{n}}
\newcommand{\R}{\mathbb{R}}
\newcommand{\N}{\mathbb{N}}

\newcommand{\p}{\partial}
\newcommand{\pn}{\partial_{\nu}}
\newcommand{\Ci}{\mathscr{C}^{\infty}}
\newcommand{\Cc}{\mathscr{C}_{c}^{\infty}}

\newenvironment{pr}{\vspace{5pt}\textbf{{\small Proof :}}\\}{\hspace{\stretch{1}}\rule{1ex}{1ex}\vspace{5pt}}
\newtheorem{thm}{Theorem}[section]

\newtheorem{lem}{Lemma}[section]

\usepackage{fancyhdr}
\fancypagestyle{plain}{\fancyhead{}}
\title{Asymptotic behavior of the transmission Euler-Bernoulli plate and wave equation with a localized Kelvin-Voigt damping}
\author{{FATHI HASSINE}\\ \textit{UR Analysis and Control of PDE UR13ES64}\\ \textit{Department of Mathematics, Faculty of Sciences of Monastir}\\ \textit{University of Monastir, 5019 Monastir, Tunisia}\\ \textit{email:} \texttt{fathi.hassine@fsm.rnu.tn}}
\date{}
\begin{document}
\maketitle
\begin{center}
\abstract{
Let a fourth and a second order evolution equations be coupled via the interface by transmission conditions, and suppose that the first one is stabilized by a localized distributed feedback. What will then be the effect of such a partial stabilization on the decay of solutions at infinity? Is the behavior of the first component sufficient to stabilize the second one? The answer given in this paper is that sufficiently smooth solutions decay logarithmically at infinity even the feedback dissipation affects an arbitrarily small open subset of the interior. The method used, in this case, is based on a frequency method, and this by combining a contradiction argument with the Carleman estimates technique to carry out a special analysis for the resolvent.
}
\end{center}
\section{Introduction}
There are several mathematical models representing physical damping. The most often encountered type of damping in vibration studies are linear viscous damping and Kelvin-Voigt damping which are special cases of proportional damping. Viscous damping usually models external friction forces such as air resistance acting on the vibrating structures and is thus called "external damping", while Kelvin-Voigt damping originate from the internal friction of the material of the vibrating structures and thus called "internal damping" or "material damping".

The study of the stabilization problem for coupled systems has attracted a lot of attention in recent years e.g.~\cite{AN,AV,D,I,Ha3,LZ,RZZ,T,ZZ1,ZZ2,ZhZh}. The systems discussed in those paper involve thermoelastic systems, fluid-structure interaction systems, and coupled wave-wave, plate-plate, or plate-wave equations. But in the case of plate-wave and for the multi-dimensional space (of interest in this paper), and as far as we know, the only models which has been treated in this subject are the model of coupled Euler-Bernoulli and wave equations with indirect damping mechanisms (see~\cite{T}) and the model arising in the control of noise, coupling the damped wave equation with a damped Kirchhoff thin plate equation (see~\cite{AN}). The system that we are going to discuss in the present paper is coupling the transversal vibration of the Euler-Bernoulli beam with Kelvin-Voigt damping distributed locally on any subdomain with the elastic wave equation.

A relevant question raised about the transmission problems and problems with locally distributed damping, is the asymptotic behavior of the solutions. Does the solution goes to zero uniformly? If this is the case, what is the rate of decay?

In~\cite{LL} and recently in~\cite{Ha1} for the transmission problem case, longitudinal and transversal vibrations of a clamped elastic beam was studied as problems with locally distributed damping. It was shown, for the one-dimensional case, that when viscoelastic damping is distributed only on a subinterval in the interior of the domain, the exponential stability holds for the transversal but not for the longitudinal motion. Besides, an optimality result was shown for longitudinal case.

\begin{figure}[htbp]
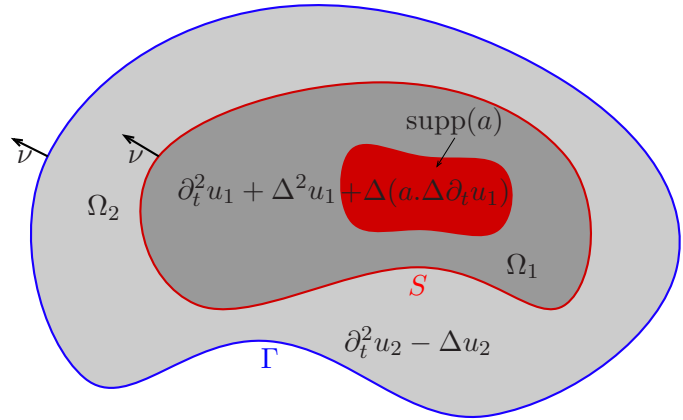

\figinit{1.4pt}
\figpt 1:(-90,30)
\figpt 2:(-30,70)\figpt 3:(50,50)
\figpt 4:(80,0)\figpt 5:(30,-40)
\figpt 6:(-30,-20)\figpt 7:(-80,-30)
\figpt 8:(-60,30)\figpt 9:(0,50)\figpt 10:(50,30)
\figpt 12:(50,-10)\figpt 13:(10,0)
\figpt 14:(-50,-10)
\figpt 15:(-8,32)\figpt 16:(12,30)\figpt 17:(33,27)
\figpt 18:(32,10)\figpt 19:(12,10)
\figpt 20:(-8,12)
\figpt 21:(20,35)\figpt 24:(15,27)
\figpt 22:(-100,35)\figpt 23:(-70,36)
\psbeginfig{}
\psset(fillmode=yes,color=0.8)
\pscurve[1,2,3,4,5,6,7,1,2,3]
\psset(color=0.6)
\pscurve[8,9,10,12,13,14,8,9,10]
\psset (color=0.8\Redrgb)
\pscurve[15,16,17,18,19,20,15,16,17]
\psset(fillmode=no,color=0.1\Bluergb)
\psset(width=0.8)
\pscurve[1,2,3,4,5,6,7,1,2,3]
\psset(fillmode=no,color=0.8\Redrgb)
\pscurve[8,9,10,12,13,14,8,9,10]
\psset(width=0.3)
\psset(fillmode=no,color=\Blackrgb)
\psset arrowhead(ratio=0.2)\psarrow[21,24]
\psset(width=0.8)
\psset arrowhead(ratio=0.3)\psarrow[1,22]\psarrow[8,23]
\psendfig
\figvisu{\figBoxA}{}{
\figwrites 16:$+\Delta\mathrm{(}a.\Delta\p_{t}u_{1}\mathrm{)}$(6)
\figwrites 13:\textcolor{red}{$S$}(2)
\figwrites 6:\textcolor{blue}{$\Gamma$}(2)
\figwritese 8:$\p_{t}^{2}u_{1}+\Delta^{2} u_{1}$(7)
\figwritese 1:$\Omega_{2}$(15)
\figwrites 13:$\p_{t}^{2}u_{2}-\Delta u_{2}$(16)
\figwriten 21:$\supp\mathrm{(}a\mathrm{)}$(0.5)
\figwritenw 12:$\Omega_{1}$(10)
\figwritew 1:$\nu$(4)
\figwritew 8:$\nu$(4)
}
\centerline{\box\figBoxA}
\caption{The domains $\Omega_{1}$ and $\Omega_{2}$ and the plate and wave operators.}
\label{fig2}
\end{figure}

Let us describe this system in detail. Let $\Omega\subset\R^{n}$ be a bounded domain with connected and smooth boundary $\Gamma$. Let $\Omega_{1}$ be a sub-domain of $\Omega$ such that $\Omega_{1}\Subset\Omega$ and set $\Omega_{2}=\Omega\backslash\overline{\Omega}_{1}$. We denote by $S$ the interface that supposed to be connected and smooth, and $\nu$ denotes the outward normal vector of $\Omega_{1}$ in $S$ and of $\Omega_{2}$ in $\Gamma$ (see Figure~\ref{fig2}). We Consider the following transmission problem
\begin{equation}\label{b1}
\left\{
\begin{array}{lll}
\p_{t}^{2}u_{1}+\Delta(\Delta u_{1}+a.\Delta\p_{t}u_{1})=0&\textrm{in}&\Omega_{1}\times]0,+\infty[,
\\
\p_{t}^{2}u_{2}-\Delta u_{2}=0&\textrm{in}&\Omega_{2}\times]0,+\infty[,
\\
u_{1}=u_{2}&\textrm{on}&S\times]0,+\infty[,
\\
\pn\Delta u_{1}+\pn u_{2}=0&\textrm{on}&S\times]0,+\infty[,
\\
\pn u_{1}=0&\textrm{on}&S\times]0,+\infty[,
\\
u_{2}=0&\textrm{on}&\Gamma\times]0,+\infty[,
\\
u_{1}(x,0)=u_{1}^{0}(x),\;\p_{t}u_{1}(x,0)=u_{1}^{1}(x)&\textrm{in}&\Omega_{1},
\\
u_{2}(x,0)=u_{2}^{0}(x),\;\p_{t}u_{2}(x,0)=u_{2}^{1}(x)&\textrm{in}&\Omega_{2}.
\end{array}\right.
\end{equation}
Where $a$ is a non negative bounded function on $\Omega_{1}$, vanishing near the interface $S$ such that there exist a non empty open domain $\omega\subset\Omega_{1}$ in such a way $a$ is strictly positive in $\overline{\omega}$.

This vibrating system is assumed to coupling the transversal and longitudinal motions (with dissipation on the plate affects an arbitrarily small open subset of its interior) through the transmission conditions as given in the third forth and fifth line of~\eqref{b1}: The first, is called the continuity condition, the second, is described by the fact that the slope of the beam is null and the third, says that the transverse force caused by the tension is equal to the transverse force due to shear. This problem was studied in~\cite{Ha2} for one-dimensional case. It was proved that the energy of the solution is decreasing with a polynomial rate for the two cases where the damping arising from the transversal motion and the damping arising from the longitudinal motion.

The energy of a solution $(u_{1},u_{2})$ of the system~\eqref{b1} at the time $t\geq 0$ is given by
\begin{equation*}
E(t)=\frac{1}{2}\left(\int_{\Omega_{1}}|\p_{t}u_{1}(x,t)|^{2}+|\Delta u_{1}(x,t)|^{2}\,\ud x+\int_{\Omega_{2}}|\p_{t}u_{2}(x,t)|^{2}+|\nabla u_{2}(x,t)|^{2}\,\ud x\right).
\end{equation*}
By means of the classical energy method, we show that
\[
E(t_{2})-E(t_{1})=-\int_{t_{1}}^{t_{2}}\!\!\!\int_{\Omega_{1}}a|\Delta\p_{t}u_{1}(x,t)|^{2}\,\ud x\,\ud t,\quad\forall\;t_{1},\,t_{2}>0.
\]
Therefore, the energy is a non-increasing function of the time variable t and our system~\eqref{b1} is dissipative. We define the Hilbert space $\mathcal{H}=X\times H$ where $H=L^{2}(\Omega_{1})\times L^{2}(\Omega_{2})$ and
\begin{equation*}
X=\Big\{(u_{1},u_{2})\in  H^{1}(\Omega_{1})\times  H_{\Gamma}^{1}(\Omega_{2}):\,\Delta u_{1}\in L^{2}(\Omega_{1}),\,u_{1\,|S}=u_{2\,|S},\,\pn u_{1\,|S}=0\Big\},
\end{equation*}
where $H_{\Gamma}^{1}(\Omega_{2})$ is the space of elements in $H^{1}(\Omega_{2})$ whose trace is zero on the boundary $\Gamma$. The space $\mathcal{H}$ is equipped with the norm
\[
\|(u_{1},u_{2},v_{1},v_{2})\|_{\mathcal{H}}^{2}=\|\Delta u_{1}\|_{L^{2}(\Omega_{1})}^{2}+\|\nabla u_{2}\|_{L^{2}(\Omega_{2})}^{2}+\|v_{1}\|_{L^{2}(\Omega_{1})}^{2}+\|v_{2}\|_{L^{2}(\Omega_{2})}^{2}.
\]
We define the operator
\[
\mathcal{A}\left(\begin{array}{c}
u_{1}
\\
u_{2}
\\
v_{1}
\\
v_{2}
\end{array}\right)=(v_{1},v_{2},-\Delta(\Delta u_{1}+a\Delta v_{1}),\Delta u_{2})
\]
whose domain is given by
\begin{align*}
\mathcal{D}(\mathcal{A})=\big\{(u_{1},u_{2},v_{1},v_{2})\in\mathcal{H}:\;(v_{1},v_{2},\Delta(\Delta u_{1}+a\Delta v_{1}),\Delta u_{2})\in\mathcal{H},\,&\pn\Delta u_{1\,|S}+\pn u_{2\,|S}=0\big\}.
\end{align*}

Our main result is the following
\begin{thm}\label{b2}
For any $k\in\N$ there exists $C>0$ such that for any initial data $(u_{0}^{0},u_{1}^{0},u_{0}^{1},u_{1}^{1})\in\mathcal{D}(\mathcal{A}^{k})$ the solution $(u_{1}(x,t),u_{2}(x,t))$ of~\eqref{b1} starting from $(u_{0}^{0},u_{1}^{0},u_{0}^{1},u_{1}^{1})$ satisfying
$$E(t)\leq\frac{C}{(\ln(2+t))^{2k}}\|(u_{0}^{0},u_{1}^{0},u_{0}^{1},u_{1}^{1})\|_{\mathcal{D}(\mathcal{A}^{k})}^{2},\quad\forall\,t>0.$$
\end{thm}
We should mention here that the subject of stabilization of transmission problems with localized Kelvin-Voigt dissipation is perhaps not intensively studied but is not new in fact, in~\cite{ARSV} the authors consider the transmission problem of a material composed by three components, one of them is a Kelvin-Voigt viscoelastic material, the second is an elastic material (no dissipation) and the third is an elastic material inserted with a frictional damping mechanism where they show different types of decay rate of energy depends on which component is in the middle, and in~\cite{BRA} the authors consider a transmission problem for the longitudinal displacement of a Euler-Bernoulli beam, where one small part of the beam is made of a viscoelastic material with Kelvin-Voigt constitutive relation in which they show that the semigroup associated to the system is exponentially stable.
\subsection*{The method of analysis}
Besides the fact that the Kelvin-Voigt damping is unbounded in the energy space and the fact that the resolvent of the system operator is not compact, the main difficulty of our problem is none other than the different speeds of propagation due to the coupling between the wave equation and the plate equation. The method that we consider here consist to the use of the Burq's result~\cite{Bur} (see also~\cite{BD}) which links, for dissipative operators, logarithmic decay to resolvent estimates with exponential loss. The main idea, as introduced by Lebeau~\cite{Leb} is to use the what's called, Carleman estimates (see also~\cite{D} for the case of non linear damping and~\cite{ET} for the case of hyperbolic systems). Unlike to the works of~\cite{B,D,I}, here Carleman estimate does not seem to be enough, that is why we have combined it with some contradiction arguments to establish the kind of resolvent estimate cited above. Moreover, to deal with the high order of the plate equation, Carleman estimate (Theorem~\ref{c4}) is established for system of second order~\eqref{c2} which is derived from the resolvent problem~\eqref{f3} by decomposing the plate equation into two second order operators~\eqref{f6}. 

The outline of this paper is as follows: In section~\ref{b3} we show that the corresponding model are well posed, in section~\ref{a1} we give the Carleman estimates and we construct a suitable weight functions that satisfy the H{\"o}rmander's assumption. In section~\ref{f1} we prove the resolvent estimate which provides the logarithmic decay result given by Theorem~\ref{b2}. 

\section{Existence and uniqueness}\label{b3}
In this section and through the semigroup theory we will show that the problem~\eqref{b1} is Well-posed. The system~\eqref{b1} can be written in the abstract form as a Cauchy problem as follows
$$\left\{\begin{array}{lcl}
\p_{t}\left(\begin{array}{c}
u_{1}
\\
u_{2}
\\
v_{1}
\\
v_{2}
\end{array}\right)(t,x)=\mathcal{A}\left(\begin{array}{c}
u_{1}
\\
u_{2}
\\
v_{1}
\\
v_{2}
\end{array}\right)(t,x)&\text{if}&(t,x)\in]0,+\infty[\times\Omega,
\\
\left(\begin{array}{l}
u_{1}
\\
u_{2}
\\
v_{1}
\\
v_{2}
\end{array}\right)(0,x)=\left(\begin{array}{l}
u_{1}^{0}
\\
u_{2}^{0}
\\
u_{1}^{1}
\\
u_{2}^{1}
\end{array}\right)(x)&\text{if}&x\in\Omega,
\end{array}\right.$$
where we recall that the operator $\mathcal{A}$ is defined by
\[
\mathcal{A}\left(\begin{array}{c}
u_{1}
\\
u_{2}
\\
v_{1}
\\
v_{2}
\end{array}\right)=(v_{1},v_{2},-\Delta(\Delta u_{1}+a\Delta v_{1}),\Delta u_{2})
\]
with domain
\begin{align*}
\mathcal{D}(\mathcal{A})=\big\{(u_{1},u_{2},v_{1},v_{2})\in\mathcal{H}:\;(v_{1},v_{2},\Delta(\Delta u_{1}+a\Delta v_{1}),\Delta u_{2})\in\mathcal{H},\,&\pn\Delta u_{1\,|S}+\pn u_{2\,|S}=0\big\}.
\end{align*}

In the space $Y=L^{2}(\Omega_{1})\times H_{\Gamma}^{1}(\Omega_{2})$ we define the operator $G$ as follows
\[
G\left(\begin{array}{l}
u_{1}
\\
u_{2}
\end{array}\right)=\left(\begin{array}{c}
-\Delta u_{1}
\\
u_{2}
\end{array}\right)
\]
with domain
\begin{equation*}
\mathcal{D}(G)=X=\Big\{(u_{1},u_{2})\in  H^{1}(\Omega_{1})\times  H_{\Gamma}^{1}(\Omega_{2}):\,\Delta u_{1}\in L^{2}(\Omega_{1}),\,u_{1\,|S}=u_{2\,|S},\,\pn u_{1\,|S}=0\Big\}.
\end{equation*}
We define a norm in the space $X$ by
$$\|(u_{1},u_{2})\|_{X}^{2}=\|\Delta u_{1}\|_{L^{2}(\Omega_{1})}^{2}+\|\nabla u_{2}\|_{L^{2}(\Omega_{2})}^{2}.$$
The graph norm of $G$ is given by
$$\|(u_{1},u_{2})\|_{gr(G)}^{2}=\|u_{1}\|_{L^{2}(\Omega_{1})}^{2}+\|\Delta u_{1}\|_{L^{2}(\Omega_{1})}^{2}+\|\nabla u_{2}\|_{L^{2}(\Omega_{2})}^{2}.$$
\begin{lem}\label{b5}
$(X,\|\,.\,\|_{X})$ is a Hilbert space with a norm equivalent to the graph norm of $G$.
\end{lem}
\begin{pr}
Let's note first, by setting $u=\mathbb{1}_{\Omega_{1}}u_{1}+\mathbb{1}_{\Omega_{2}}u_{2}$, that the continuity transmission condition $u_{1\,|S}=u_{2\,|S}$ allows to look at $u$ as an element of $H^{1}(\Omega)$. Hence by Green's formula and Poincar\'e inequality there exists $C>0$ such that for every $(u_{1},u_{2})\in X$ we have
\begin{equation*}
\left\langle G\left(\begin{array}{c}
u_{1}
\\
u_{2}
\end{array}\right),\left(\begin{array}{c}
u_{1}
\\
u_{2}
\end{array}\right)\right\rangle_{Y}=\|\nabla u_{1}\|_{L^{2}(\Omega_{1})}^{2}+\|\nabla u_{2}\|_{L^{2}(\Omega_{2})}^{2}\geq C(\|u_{1}\|_{L^{2}(\Omega_{1})}^{2}+\|\nabla u_{2}\|_{L^{2}(\Omega_{2})}^{2}).
\end{equation*}
In particular $G$ is a strictly positive operator on $Y$. Besides, since for every $(u_{1},u_{2})\in X$ we have
$$\|(u_{1},u_{2})\|_{X}.\|(u_{1},u_{2})\|_{Y}=\|G(u_{1},u_{2})\|_{Y}.\|(u_{1},u_{2})\|_{Y}\geq\left\langle G\left(\begin{array}{c}
u_{1}
\\
u_{2}
\end{array}\right),\left(\begin{array}{c}
u_{1}
\\
u_{2}
\end{array}\right)\right\rangle_{Y}\geq C\|(u_{1},u_{2})\|_{Y}^{2},$$
then the equivalence of the two norms holds.
\\
To finish the proof we have only to prove that $G$ is a closed operator on $Y$. Let $(u_{1,n},u_{2,n})\in Y$, $(u_{1},u_{2})$ and $(f_{1},f_{2})$ such that
$$
u_{1,n}\longrightarrow u_{1},\; -\Delta u_{1,n}\longrightarrow f_{1}\;\text{in}\;L^{2}(\Omega_{1})\;\text{and}\;u_{2,n}\longrightarrow f_{2},\; u_{2,n}\longrightarrow u_{2}\;\text{in}\;H_{\Gamma}^{1}(\Omega_{2})\;\text{as}\;n\longrightarrow+\infty.
$$
Therefore, $u_{2}=f_{2}\in H_{\Gamma}^{1}(\Omega_{2})$ and in since for all $\vfi\in \Cc(\Omega_{1})$,
$$
\langle-f_{1},\vfi\rangle_{D'(\Omega_{1})}=\!\!\lim_{n\to+\infty}\langle\Delta u_{1,n},\vfi\rangle_{D'(\Omega_{1})}=\!\!\lim_{n\to+\infty}\langle u_{1,n},\Delta\vfi\rangle_{D'(\Omega_{1})}=\langle u_{1},\Delta\vfi\rangle_{D'(\Omega_{1})}=\langle\Delta u_{1},\vfi\rangle_{D'(\Omega_{1})},
$$
then we obtain also $-\Delta u_{1}=f_{1}\in L^{2}(\Omega_{1})$. In the other hand, while
$$
\|\nabla(u_{1,n}-u_{1,m})\|_{L^{2}(\Omega_{1})}=-\langle\Delta(u_{1,n}-u_{1,m}),u_{1,n}-u_{1,m}\rangle_{L^{2}(\Omega_{1})}^{2}\longrightarrow0\;\text{as}\; n,m\longrightarrow+\infty
$$
then $u_{1,n}$ is a Cauchy sequence in $H^{1}(\Omega_{1})$, that converge to $u_{1}$ in $H^{1}(\Omega_{1})$ where we argue this fact as follows,
$$
\lim_{n\to+\infty}\langle\nabla u_{1,n},\vfi\rangle_{D'(\Omega_{1})}=-\!\!\lim_{n\to+\infty}\langle u_{1,n},\nabla\vfi\rangle_{D'(\Omega_{1})}=-\langle u_{1},\nabla\vfi\rangle_{D'(\Omega_{1})}=\langle\nabla u_{1},\vfi\rangle_{D'(\Omega_{1})}.
$$
For the transmission conditions we have
$$
\|u_{1}-u_{2}\|_{H^{\frac{1}{2}}(S)}=\|u_{1}-u_{1,n}+u_{2,n}-u_{2}\|_{H^{\frac{1}{2}}(S)}\leq C(\|u_{1}-u_{1,n}\|_{H^{1}(\Omega_{1})}+\|u_{2}-u_{2,n}\|_{H^{1}(\Omega_{2})})\longrightarrow0.
$$
and
$$
\|\pn u_{1}\|_{H^{-\frac{1}{2}}(S)}=\|\pn(u_{1}-u_{1,n})\|_{H^{-\frac{1}{2}}(S)}\leq C(\|u_{1}-u_{1,n}\|_{H^{1}(\Omega_{1})}+\|\Delta u_{1}-\Delta u_{1,n}\|_{L^{2}(\Omega_{1})})\longrightarrow0,
$$
where we have used here~\cite[Theorem 13.7.6]{TW}. This show now that $G$ is a closed operator and this conclude the proof.
\end{pr}
\begin{thm}
The operator $\mathcal{A}$ is m-dissipative and especially it generates a strongly semigroup of contraction in $\mathcal{H}$.
\end{thm}
\begin{pr}
According to Lumer-Phillips theorem (see~\cite[Theorem 3.8.4]{TW}) we have only to prove that $\mathcal{A}$ is m-dissipative.
\\
Let $(u_{1},u_{2},v_{1},v_{2})\in\mathcal{D}(\mathcal{A})$ then by Green's formula we have
\begin{eqnarray*}
\re\left\langle\mathcal{A}\left(\begin{array}{c}
u_{1}
\\
u_{2}
\\
v_{1}
\\
v_{2}
\end{array}\right),\left(\begin{array}{c}
u_{1}
\\
u_{2}
\\
v_{1}
\\
v_{2}
\end{array}\right)\right\rangle_{\mathcal{H}}&=&\re\left\langle\left(\begin{array}{c}
v_{1}
\\
v_{2}
\\
-\Delta(\Delta u_{1}+a\Delta v_{1})
\\
\Delta u_{2}
\end{array}\right),\left(\begin{array}{l}
u_{1}
\\
u_{2}
\\
v_{1}
\\
v_{2}
\end{array}\right)\right\rangle_{\mathcal{H}}
\\
&=&-\|a^{\frac{1}{2}}\Delta v_{1}\|_{L^{2}(\Omega_{1})}^{2}\leq 0.
\end{eqnarray*}
This shows that $\mathcal{A}$ is dissipative.
\\
Let $(f_{1},f_{2},g_{1},g_{2})\in\mathcal{H}$ and let's find a quadruplet $(u_{1},u_{2},v_{1},v_{2})\in\mathcal{D}(\mathcal{A})$ such that
$$\left(\id-\mathcal{A}\right)\left(\begin{array}{c}
u_{1}
\\
u_{2}
\\
v_{1}
\\
v_{2}
\end{array}\right)=\left(\begin{array}{c}
u_{1}-v_{1}
\\
u_{2}-v_{2}
\\
v_{1}+\Delta(\Delta u_{1}+a\Delta v_{1})
\\
v_{2}-\Delta u_{2}
\end{array}\right)=\left(\begin{array}{c}
f_{1}
\\
f_{2}
\\
g_{1}
\\
g_{2}
\end{array}\right),$$
This amounts to finding $(u_{1},u_{2},v_{1},v_{2})\in\mathcal{D}(\mathcal{A})$ that satisfies the following system
\begin{equation*}
\left\{\begin{array}{lll}
v_{1}=u_{1}-f_{1}&\text{in}&\Omega_{1}
\\
v_{2}=u_{2}-f_{2}&\text{in}&\Omega_{2}
\\
u_{1}+\Delta((1+a)\Delta u_{1}-a\Delta f_{1})=f_{1}+g_{1}&\text{in}&\Omega_{1}
\\
u_{2}-\Delta u_{2}=f_{2}+g_{2}&\text{in}&\Omega_{2}.
\end{array}\right.
\end{equation*}
From Lemma~\ref{b5} and the Riesz representation theorem, we can find a unique $(u_{1},u_{2})\in X$ such that for all $(\vfi_{1},\vfi_{2})\in X$ we have
\begin{equation*}
\begin{split}
\langle f_{1}+g_{1},\vfi_{1}\rangle_{L^{2}(\Omega_{1})}+\langle f_{2}+g_{2},\vfi_{2}\rangle_{L^{2}(\Omega_{2})}+\langle a\Delta f_{1},\Delta\vfi_{1}\rangle_{L^{2}(\Omega_{1})}=\langle u_{1},\vfi_{1}\rangle_{L^{2}(\Omega_{1})}
\\
+\langle(1+a)\Delta u_{1},\Delta\vfi_{1}\rangle_{L^{2}(\Omega_{1})}+\langle u_{2},\vfi_{2}\rangle_{L^{2}(\Omega_{2})}+\langle\nabla u_{2},\nabla\vfi_{2}\rangle_{L^{2}(\Omega_{2})}.
\end{split}
\end{equation*}
Then by Green's formula we obtain
\begin{equation}\label{b7}
\begin{split}
\langle\Delta((1+a)\Delta u_{1}-a\Delta f_{1})+(u_{1}-f_{1}-g_{1}),\vfi_{1}\rangle_{L^{2}(\Omega_{1})}+\langle f_{2}+g_{2}+\Delta u_{2}-u_{2},\vfi_{2}\rangle_{L^{2}(\Omega_{2})}
\\
=\langle\pn u_{2}+\pn\Delta u_{1},\vfi_{1}\rangle_{L^{2}(S)}.
\end{split}
\end{equation}
In particular for all $(\vfi_{1},\vfi_{2})\in \Cc(\Omega_{1})\times\Cc(\Omega_{2})$ we have
$$\langle\Delta((1+a)\Delta u_{1}-a\Delta f_{1})+(u_{1}-f_{1}-g_{1}),\vfi_{1}\rangle_{L^{2}(\Omega_{1})}+\langle f_{2}+g_{2}+\Delta u_{2}-u_{2},\vfi_{2}\rangle_{L^{2}(\Omega_{2})}=0$$
then we find
\begin{equation}\label{b8}
\left\{\begin{array}{ll}
u_{1}+\Delta((1+a)\Delta u_{1}-a\Delta f_{1})=f_{1}+g_{1}&\text{in}\; L^{2}(\Omega_{1}),
\\
u_{2}-\Delta u_{2}=f_{2}+g_{2}&\text{in}\; L^{2}(\Omega_{2}).
\end{array}\right.
\end{equation}
Then from~\eqref{b7} and~\eqref{b8} we obtain
\begin{align*}
\langle\pn u_{2}+\pn\Delta u_{1},\vfi_{1}\rangle_{L^{2}(S)}=0
\end{align*}
and this show the following equality
$$\pn \Delta u_{1\,|S}+\pn u_{2\,|S}=0.$$
And this give end to our proof.
\end{pr}

One consequence of this last result is that if we assume that $(u_{1}^{0},u_{2}^{0},u_{1}^{1},u_{2}^{1})\in\mathcal{D}(\mathcal{A})$, there exists a unique solution of~\eqref{b1} which can be expressed by means of a semigroup on $\mathcal{H}$ as follows
\begin{equation}\label{b9}
\left(\begin{array}{c}
u_{1}
\\
u_{2}
\\
\p_{t}u_{1}
\\
\p_{t}u_{2}
\end{array}\right)=e^{t\mathcal{A}}\left(\begin{array}{c}
u_{1}^{0}
\\
u_{2}^{0}
\\
u_{1}^{1}
\\
u_{2}^{1}
\end{array}\right)
\end{equation}
where $e^{t\mathcal{A}}$ is the semigroup of the operator $\mathcal{A}$. And we have the following regularity of the solution
$$
\left(\begin{array}{c}
u_{1}
\\
u_{2}
\\
\p_{t}u_{1}
\\
\p_{t}u_{2}\end{array}\right)
\in C([0,+\infty[,\mathcal{D}(\mathcal{A}))\cap C^{1}([0,+\infty[,\mathcal{H}).$$
\\
And if $(u_{1}^{0},u_{2}^{0},u_{1}^{1},u_{2}^{1})\in\mathcal{H}$, the function $(u_{1}(t),u_{2}(t))$ given by~\eqref{b9} is the mild solution of~\eqref{b1}.
\section{Carleman estimate near the surface}\label{a1}
This section is devoted to establish the Carleman estimate.

We set the operator
$$A(x,D):=\left\{\begin{array}{ll}
\displaystyle A_{1}(x,D):=-\p_{x_{n}}^{2}+R_{1}(x,\p_{x'}/i)\pm\frac{1}{h}& x_{n}>0
\\
\\
A_{2}(x,D):=-\p_{x_{n}}^{2}+R_{2}(x,\p_{x'}/i)-\displaystyle\frac{1}{h^{2}}& x_{n}<0
\end{array}\right.$$
with $h$ is a small semi-classical parameter and where
$$R(x,\xi')=\left\{\begin{array}{ll}
R_{1}(x,\xi')& x_{n}>0
\\
R_{2}(x,\xi')& x_{n}<0
\end{array}\right.$$
is a second order polynomial in $\xi'$ with coefficients in $\R$ with principal symbol
$$r(x,\xi')=\left\{\begin{array}{ll}
r_{1}(x,\xi')& x_{n}>0
\\
r_{2}(x,\xi')& x_{n}<0
\end{array}\right.$$
that satisfy
\begin{equation*}
\left\{\begin{array}{l}
r_{1}(x,\xi')\geq C|\xi'|^{2}\quad\forall\,x_{n}>0,\;\forall\,\xi'\in\R^{n-1}
\\
r_{2}(x,\xi')\geq C|\xi'|^{2}\quad\forall\,x_{n}<0,\;\forall\,\xi'\in\R^{n-1}.
\end{array}\right.
\end{equation*}

We consider the following transmission problem
\begin{equation}\label{a3}
\left\{\begin{array}{ll}
A(x,D)(w)=f& x_{n}\neq0
\\
hw(x',0^{+})= w(x',0^{-})+\theta&
\\
\pn w(x',0^{+})=\pn w(x',0^{+})+\Theta.&
\end{array}\right.
\end{equation}

Let $V=V'\times]-\epsilon,\epsilon[$ be an open set of $\R^{n}$, follow to~\cite{RR1} we set 
$$\R_{-}^{n}=\{x:\;x_{n}<0\},\quad\R_{+}^{n}=\{x:\;x_{n}>0\},\quad V^{g}=V\cap\R_{-}^{n},\quad V^{d}=V\cap\R_{+}^{n}.$$
For a compact set $K$ of $V$ we set $K^{g}=K\cap\overline{\R_{-}^{n}}$ and $K^{d}=K\cap\overline{\R_{+}^{n}}$. We then denote by $\Cc(K^{d})$ (resp. $\Cc(K^{g})$) the space of functions that are $\Ci$ in $\overline{\R_{+}^{n}}$ (resp. $\overline{\R_{-}^{n}}$) with support in $K^{d}$ (resp. $K^{g}$). 

We let $\vfi$ a weight function and we define in both side of $S$ the conjugate operator
$$A_{\vfi}=h^{2}\e^{\vfi/h}A\e^{-\vfi/h}$$
with principal symbol 
$$
a_{\vfi}(x,\xi)=\left\{\begin{array}{ll}
(\xi_{n}+i(\p_{x_{n}}\vfi))^{2}+r_{1}(x,\xi'+i(\p_{x'}\vfi))& x_{n}>0
\\
(\xi_{n}+i(\p_{x_{n}}\vfi))^{2}+r_{2}(x,\xi'+i(\p_{x'}\vfi))-1& x_{n}<0.
\end{array}\right.
$$
We suppose that the weight function $\vfi$ is in $\mathscr{C}(\overline{V})$, $\vfi_{|\overline{\Rn_{-}}}\in\Ci(\overline{V}^{g})$, $\vfi_{|\overline{\Rn_{+}}}\in\Ci(\overline{V}^{d})$ and such that
\begin{enumerate}
	\item $|\nabla\varphi|(x)>0$ in $\overline{V}$.
	\item For all $x'\in V'$
	\begin{equation}\label{a2}
	\left\{
	\begin{array}{l}
	(\p_{x_{n}}\vfi)(x',0^{+})>0 \textrm{ and }(\p_{x_{n}}\vfi)(x',0^{-})>0
	\\
	(\p_{x_{n}}\vfi)^{2}(x',0^{+})-(\p_{x_{n}}\vfi)^{2}(x',0^{-})>1.
	\end{array}\right.
	\end{equation}
	\item The sub-ellipticity condition:
	\begin{equation}\label{a5}
	\forall\,(x,\xi)\in\overline{V}\times\Rn;\; a_{\vfi}(x,\xi)=0\Longrightarrow\{\re(a_{\vfi}),\im(a_{\vfi})\}(x,\xi)>0.
	\end{equation}
\end{enumerate}
Follows to Le Rousseau and Robbiano result~\cite{RR1} we can prove by using the same argument and the exactly the same steps we can prove the following result
\begin{thm}\label{a4}
Let $K$ be a compact subset of $V$ and $\vfi$ a weight function satisfying the above assumption, then there exist $C>0$ and $h_{0}>0$ such that
\begin{equation*}
\begin{split}
h\|\e^{\vfi/h}w\|_{0}^{2}+h^{3}\|\e^{\vfi/h}\nabla w\|_{0}^{2}&+h|\e^{\vfi/h}w_{|x_{n}=0^{\pm}}|_{0}^{2}+h^{3}|\e^{\vfi/h}\nabla_{x'} w_{|x_{n}=0^{\pm}}|_{0}^{2}+h^{3}|\e^{\vfi/h}\p_{x_{n}}w_{|x_{n}=0^{\pm}}|_{0}^{2}
\\
&\leq C(h^{4}\|\e^{\vfi/h}f\|_{0}^{2}+h|\e^{\vfi/h}\theta|_{0}^{2}+h^{3}|\e^{\vfi/h}\nabla_{x'}\theta|_{0}^{2}+h^{3}|\e^{\vfi/h}\Theta|_{0}^{2})
\end{split}
\end{equation*}
for all $0<h\leq h_{0}$, $w$ and $f$ satisfying~\eqref{a3} where $w_{|\overline{\Rn_{-}}}\in\Cc(K^{g})$ and $w_{|\overline{\Rn_{+}}}\in\Cc(K^{d})$.
\end{thm}

The proof of Carleman estimate is the same for both, in this paper and in the Le Rousseau and Robbiano paper~\cite{RR1} even that the transmission conditions are different. In fact, while in~\cite{RR1} they are depending on some diffusion coefficients where an additional assumption (\cite[(2.2)]{RR1}), on the jump at the interface of the weight functions, is assumed in addition to that given above, here the transmission conditions depend on the pseudo-differential parameter $h$ where, for $h$ small enough this scaling coefficient allows us to ensure the assumption of Le Rousseau and Robbiano~\cite[(2.2)]{RR1} which became in our case $(\p_{x_{n}}\vfi)(x',0^{+})-h(\p_{x_{n}}\vfi)(x',0^{-})>0$. Thus we may notice how the scaling coefficient $h$ is playing the same role as the diffusion coefficients in the Rousseau and Robbiano paper~\cite{RR1}. Noting also another version of this analysis appeared more recently in~\cite{RLR}.

The purpose of the rest of this section is to provide a global Carleman estimate for a transmission problem with three entries governed by a three elliptic operators. Besides, we will try to construct a suitable weight functions that will be needed in the next section.

\begin{figure}[htbp]
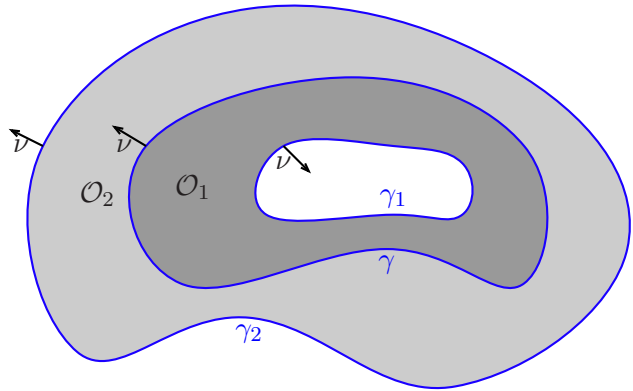

\figinit{1.3pt}
\figpt 1:(-90,30)
\figpt 2:(-30,70)\figpt 3:(50,50)
\figpt 4:(80,0)\figpt 5:(30,-40)
\figpt 6:(-30,-20)\figpt 7:(-80,-30)
\figpt 8:(-60,30)\figpt 9:(0,50)\figpt 10:(50,30)
\figpt 12:(50,-10)\figpt 13:(10,0)
\figpt 14:(-50,-10)
\figpt 15:(-20,30)\figpt 16:(12,30)\figpt 17:(32,25)
\figpt 18:(32,10)\figpt 19:(12,10)
\figpt 20:(-25,10)
\figpt 21:(-100,35)\figpt 22:(-70,36)\figpt 23:(-12,22)
\psbeginfig{}
\psset(fillmode=yes,color=0.8)
\pscurve[1,2,3,4,5,6,7,1,2,3]
\psset(color=0.6)
\pscurve[8,9,10,12,13,14,8,9,10]
\psset (color=1)
\pscurve[15,16,17,18,19,20,15,16,17]
\psset(fillmode=no,color=0.1\Bluergb)
\psset(width=0.8)
\pscurve[1,2,3,4,5,6,7,1,2,3]
\pscurve[8,9,10,12,13,14,8,9,10]
\pscurve[15,16,17,18,19,20,15,16,17]
\psset(fillmode=no,color=\Blackrgb)
\psset arrowhead(ratio=0.3)\psarrow[1,21]\psarrow[8,22]\psarrow[15,23]
\psendfig
\figvisu{\figBoxB}{}{
\figwrites 13:\textcolor{blue}{$\gamma$}(2)
\figwrites 6:\textcolor{blue}{$\gamma_{2}$}(2)
\figwriten 19:\textcolor{blue}{$\gamma_{1}$}(2)
\figwritese 1:$\mathcal{O}_{2}$(15)
\figwritese 8:$\mathcal{O}_{1}$(12)
\figwritew 1:$\nu$(4)
\figwritew 8:$\nu$(4)
\figwrites 15:$\nu$(3)
}
\centerline{\box\figBoxB}
\caption{The domains $\mathcal{O}_{1}$ and $\mathcal{O}_{2}$.}
\label{fig3}
\end{figure}

Let $\mathcal{O}_{1}$ and $\mathcal{O}_{2}$ be two open and disjoint domains with smooth boundary and we suppose that $\p\mathcal{O}_{1}=\gamma\cup\gamma_{1}$ and $\p\mathcal{O}_{2}=\gamma\cup\gamma_{2}$ such that $\overline{\gamma}_{1}\cap\overline{\gamma}=\overline{\gamma}_{2}\cap\overline{\gamma}=\overline{\gamma}_{1}\cap\overline{\gamma}_{2}=\emptyset$. We denote by $\nu$ the outward normal vector of $\mathcal{O}_{1}$ in $\gamma$ and $\gamma_{1}$ and of $\Omega_{2}$ in $\gamma_{2}$ (see Figure~\ref{fig3}). We consider the following boundary and transmission value problem
\begin{equation}\label{c2}
\left\{\begin{array}{lll}
\displaystyle A_{1}'y_{1}'=f_{1}'&\text{in}&\mathcal{O}_{1}
\\
\displaystyle A_{1}''y_{1}''=f_{1}''&\text{in}&\mathcal{O}_{1}
\\
\displaystyle A_{2}y_{2}=f_{2}&\text{in}&\mathcal{O}_{2}
\\
hy_{1}'=y_{2}+hy_{1}''&\text{on}&\gamma
\\
\pn y_{1}'=\pn y_{2}&\text{on}&\gamma
\\
hy_{1}''=y_{2}+\theta&\text{on}&\gamma
\\
\pn y_{1}''=\pn y_{2}&\text{on}&\gamma
\\
y_{2}=0&\text{on}&\gamma_{2}.
\end{array}\right.
\end{equation}
where $A_{1}'$, $A_{2}''$ and $A_{2}$ are differential operators defined by
$$
\displaystyle A_{1}'y_{1}'=-\Delta y_{1}'-y_{1}'/h,\quad A_{1}''y_{1}''=-\Delta y_{1}''+y_{1}''/h,\quad A_{2}y_{2}=-\Delta y_{2}-y_{2}/h^{2}.
$$

We define the conjugate operators of $A_{1}'$, $A_{1}''$ and $A_{2}$ respectively by
\begin{equation*}
\begin{array}{c}
A_{\vfi_{1}}'(x,D,h)=h^{2}\e^{\vfi_{1}/h}A_{1}'\e^{-\vfi_{1}/h},\;
A_{\vfi_{1}}''(x,D,h)=h^{2}\e^{\vfi_{1}/h}A_{1}''\e^{-\vfi_{1}/h},
\\
A_{\vfi_{2}}(x,D,h)=h^{2}\e^{\vfi_{2}/h}A_{2}\e^{-\vfi_{2}/h},
\end{array}
\end{equation*}
where $A_{\vfi_{1}}'(x,D,h)$ and $A_{\vfi_{1}}''(x,D,h)$ are of principal symbol
$$
a_{\vfi_{1}}(x,\xi)=|\xi+i\nabla\vfi_{1}|^{2}
$$
and that of $A_{\vfi_{2}}(x,D,h)$ is
$$
a_{\vfi_{2}}(x,\xi)=|\xi+i\nabla\vfi_{2}|^{2}-1
$$
where $\vfi_{1}$ and $\vfi_{2}$ are two weight functions defined respectively in $\overline{\mathcal{O}}_{1}$ and $\overline{\mathcal{O}}_{2}$.

We suppose that the weight functions $\vfi_{1}\in\Ci(\overline{\mathcal{O}}_{1})$, $\vfi_{2}\in\Ci(\overline{\mathcal{O}}_{2})$ and satisfy
\begin{enumerate}
	\item $\vfi_{1\,|\gamma}=\vfi_{2\,|\gamma}$.
	\item $|\nabla\vfi_{1}|(x)>0$ in $\overline{\mathcal{O}}_{1}$ and $|\nabla\vfi_{2}|(x)>0$ in $\overline{\mathcal{O}}_{2}$.
	\item $\pn\vfi_{1\,|\gamma}<0$ and $\pn\vfi_{2\,|\gamma}<0$.
	\item $(\pn\vfi_{1\,|\gamma})^{2}-(\pn\vfi_{2\,|\gamma})^{2}>1$
	\item $\pn\vfi_{1\,|\gamma_{1}}\neq0$ and $\pn\vfi_{2\,|\gamma_{2}}<0$.
	\item The sub-ellipticity condition:
	\begin{equation}\label{c5}
	\left\{\begin{array}{l}
	\forall\,(x,\xi)\in\overline{\mathcal{O}}_{1}\times\Rn;\; a_{\vfi_{1}}(x,\xi)=0\Longrightarrow\{\re(a_{\vfi_{1}}),\im(a_{\vfi_{1}})\}(x,\xi)>0
	\\
	\forall\,(x,\xi)\in\overline{\mathcal{O}}_{2}\times\Rn;\; a_{\vfi_{2}}(x,\xi)=0\Longrightarrow\{\re(a_{\vfi_{2}}),\im(a_{\vfi_{2}})\}(x,\xi)>0.
	\end{array}\right.
	\end{equation}
\end{enumerate}
Under these assumption the global Carleman estimate is given by the following
\begin{thm}\label{c4}
Let $\vfi_{1}$ and $\vfi_{2}$ the two weight functions as described above, then there exist $C>0$ and $h_{0}>0$ such that
\begin{equation*}
\begin{split}
h\|\e^{\vfi_{1}/h}y_{1}'\|_{L^{2}(\mathcal{O}_{1})}^{2}+h^{3}\|\e^{\vfi_{1}/h}\nabla y_{1}'\|_{L^{2}(\mathcal{O}_{1})}^{2}+h|\e^{\vfi_{1}/h}y_{1}'|_{L^{2}(\gamma)}^{2}+h^{3}|\e^{\vfi_{1}/h}\nabla y_{1}'|_{L^{2}(\gamma)}^{2}
\\
+h^{3}|\e^{\vfi_{1}/h}\pn y_{1}'|_{L^{2}(\gamma)}^{2}+h\|\e^{\vfi_{1}/h}y_{1}''\|_{L^{2}(\mathcal{O}_{1})}^{2}
+h^{3}\|\e^{\vfi_{1}/h}\nabla y_{1}''\|_{L^{2}(\mathcal{O}_{1})}^{2}+h|\e^{\vfi_{1}/h}y_{1}''|_{L^{2}(\gamma)}^{2}
\\
+h^{3}|\e^{\vfi_{1}/h}\nabla y_{1}''|_{L^{2}(\gamma)}^{2}+h^{3}|\e^{\vfi_{1}/h}\pn y_{1}''|_{L^{2}(\gamma)}^{2}+h\|\e^{\vfi_{2}/h}y_{2}\|_{L^{2}(\mathcal{O}_{2})}^{2}
+h^{3}\|\e^{\vfi_{2}/h}\nabla y_{2}\|_{L^{2}(\mathcal{O}_{2})}^{2}
\\
+h|\e^{\vfi_{2}/h}y_{2}|_{L^{2}(\gamma)}^{2}+h^{3}|\e^{\vfi_{2}/h}\nabla y_{2}|_{L^{2}(\gamma)}^{2}+h^{3}|\e^{\vfi_{2}/h}\pn y_{2}|_{L^{2}(\gamma)}^{2}\leq C(h^{4}\|\e^{\vfi_{1}/h}f_{1}'\|_{L^{2}(\mathcal{O}_{1})}^{2}
\\
+h^{4}\|\e^{\vfi_{1}/h}f_{1}''\|_{L^{2}(\mathcal{O}_{1})}^{2}+h^{4}\|\e^{\vfi_{1}/h}f_{2}\|_{L^{2}(\mathcal{O}_{2})}^{2}+h|\e^{\vfi/h}\theta|_{L^{2}(\gamma)}^{2}+h^{3}|\e^{\vfi/h}\nabla\theta|_{L^{2}(\gamma)}^{2}
\\
+h|\e^{\vfi_{1}/h}y_{1}'|_{L^{2}(\gamma_{1})}^{2}+h^{3}|\e^{\vfi_{1}/h}\pn y_{1}'|_{L^{2}(\gamma_{1})}^{2}+h|\e^{\vfi_{1}/h}y_{1}''|_{L^{2}(\gamma_{1})}^{2}+h^{3}|\e^{\vfi_{1}/h}\pn y_{1}''|_{L^{2}(\gamma_{1})}^{2})
\end{split}
\end{equation*}
for all $0<h\leq h_{0}$ and $y_{1}',\,y_{1}''\in\Ci(\overline{\mathcal{O}}_{1})$, $y_{2}\in\Ci(\overline{\mathcal{O}}_{2})$, $f_{1}',\,f_{1}''$ and $f_{2}$ satisfying to the system~\eqref{c2}.
\end{thm}
\begin{pr}
The proof follows easily from Theorem~\ref{a4} in fact, system~\eqref{c2} can be shown as a combination of two transmission problems, the first is by consider only the equation with entries $y_{1}'$ and $y_{2}$ only, where in the first transmission equation the term $\theta'=hy_{1}''$ should be seen as an error in the continuity of the trace of $y_{1}'$ and $y_{2}$, namely we have
\begin{equation*}
\left\{\begin{array}{lll}
\displaystyle A_{1}'y_{1}'=f_{1}'&\text{in}&\mathcal{O}_{1}
\\
\displaystyle A_{2}y_{2}=f_{2}&\text{in}&\mathcal{O}_{2}
\\
hy_{1}'=y_{2}+\theta'&\text{on}&\gamma
\\
\pn y_{1}'=\pn y_{2}&\text{on}&\gamma
\\
\pn y_{1}''=\pn y_{2}&\text{on}&\gamma
\\
y_{2}=0&\text{on}&\gamma_{2}.
\end{array}\right.
\end{equation*}
and the second problem is  with the entries $y_{1}''$ and $y_{2}$ only as follow
\begin{equation*}
\left\{\begin{array}{lll}
\displaystyle A_{1}''y_{1}''=f_{1}''&\text{in}&\mathcal{O}_{1}
\\
\displaystyle A_{2}y_{2}=f_{2}&\text{in}&\mathcal{O}_{2}
\\
hy_{1}''=y_{2}+\theta&\text{on}&\gamma
\\
\pn y_{1}''=\pn y_{2}&\text{on}&\gamma
\\
y_{2}=0&\text{on}&\gamma_{2}.
\end{array}\right.
\end{equation*}
We apply Theorem~\ref{a4} for each of these systems by taking into account \cite[Proposition 1]{LeRo1} and \cite[Proposition 2]{LeRo2} then we get
\begin{equation*}
\begin{split}
h\|\e^{\vfi_{1}/h}y_{1}'\|_{L^{2}(\mathcal{O}_{1})}^{2}+h^{3}\|\e^{\vfi_{1}/h}\nabla y_{1}'\|_{L^{2}(\mathcal{O}_{1})}^{2}+h|\e^{\vfi_{1}/h}y_{1}'|_{L^{2}(\gamma)}^{2}+h^{3}|\e^{\vfi_{1}/h}\nabla y_{1}'|_{L^{2}(\gamma)}^{2}
\\
+h^{3}|\e^{\vfi_{1}/h}\pn y_{1}'|_{L^{2}(\gamma)}^{2}+h\|\e^{\vfi_{2}/h}y_{2}\|_{L^{2}(\mathcal{O}_{2})}^{2}
+h^{3}\|\e^{\vfi_{2}/h}\nabla y_{2}\|_{L^{2}(\mathcal{O}_{2})}^{2}+h|\e^{\vfi_{2}/h}y_{2}|_{L^{2}(\gamma)}^{2}
\\
+h^{3}|\e^{\vfi_{2}/h}\nabla y_{2}|_{L^{2}(\gamma)}^{2}+h^{3}|\e^{\vfi_{2}/h}\pn y_{2}|_{L^{2}(\gamma)}^{2}\leq C(h^{4}\|\e^{\vfi_{1}/h}f_{1}'\|_{L^{2}(\mathcal{O}_{1})}^{2}+h^{4}\|\e^{\vfi_{1}/h}f_{2}\|_{L^{2}(\mathcal{O}_{2})}^{2}
\\
+h|\e^{\vfi/h}\theta'|_{L^{2}(\gamma)}^{2}+h^{3}|\e^{\vfi/h}\nabla\theta'|_{L^{2}(\gamma)}^{2}+h|\e^{\vfi_{1}/h}y_{1}'|_{L^{2}(\gamma_{1})}^{2}+h^{3}|\e^{\vfi_{1}/h}\pn y_{1}'|_{L^{2}(\gamma_{1})}^{2}),
\end{split}
\end{equation*}
and
\begin{equation*}
\begin{split}
h\|\e^{\vfi_{1}/h}y_{1}''\|_{L^{2}(\mathcal{O}_{1})}^{2}+h^{3}\|\e^{\vfi_{1}/h}\nabla y_{1}''\|_{L^{2}(\mathcal{O}_{1})}^{2}+h|\e^{\vfi_{1}/h}y_{1}''|_{L^{2}(\gamma)}^{2}+h^{3}|\e^{\vfi_{1}/h}\nabla y_{1}''|_{L^{2}(\gamma)}^{2}
\\
+h^{3}|\e^{\vfi_{1}/h}\pn y_{1}''|_{L^{2}(\gamma)}^{2}+h\|\e^{\vfi_{2}/h}y_{2}\|_{L^{2}(\mathcal{O}_{2})}^{2}
+h^{3}\|\e^{\vfi_{2}/h}\nabla y_{2}\|_{L^{2}(\mathcal{O}_{2})}^{2}+h|\e^{\vfi_{2}/h}y_{2}|_{L^{2}(\gamma)}^{2}
\\
+h^{3}|\e^{\vfi_{2}/h}\nabla y_{2}|_{L^{2}(\gamma)}^{2}+h^{3}|\e^{\vfi_{2}/h}\pn y_{2}|_{L^{2}(\gamma)}^{2}\leq C(h^{4}\|\e^{\vfi_{1}/h}f_{1}''\|_{L^{2}(\mathcal{O}_{1})}^{2}+h^{4}\|\e^{\vfi_{1}/h}f_{2}\|_{L^{2}(\mathcal{O}_{2})}^{2}
\\
+h|\e^{\vfi/h}\theta|_{L^{2}(\gamma)}^{2}+h^{3}|\e^{\vfi/h}\nabla\theta|_{L^{2}(\gamma)}^{2}+h|\e^{\vfi_{1}/h}y_{1}''|_{L^{2}(\gamma_{1})}^{2}+h^{3}|\e^{\vfi_{1}/h}\pn y_{1}''|_{L^{2}(\gamma_{1})}^{2}).
\end{split}
\end{equation*}
The result follows easily now by combing the two last estimates where the terms $\theta'=hy_{1}''$ are absorbed to the left hand side for $h>0$ small enough. 
\end{pr}

We return now to our geometric baseline as described in the introduction of this paper and we denote by $\widetilde{\Omega}_{1}=\Omega_{1}\backslash \overline{B}_{r}$ where $B_{r}$ is an open ball of $\Omega_{1}$ with radius $r>0$ such that $\overline{B}_{r}\subset\Omega_{1}$. We try to find four phases $\foo,\,\phot,\,\phto$ and $\ftt$ satisfying the H{\"o}rmander's condition except in a finite number of balls where $\foo$ or $\phot$ (resp. $\phto$ or $\ftt$) do not satisfy this condition the other does and is strictly greater. 

According to~\cite[Proposition 3.2]{Bur} we can find two $\Ci$ functions $\foo$ and $\phot$ in $\overline{\widetilde{\Omega}}_{1}$ (resp. $\phto$ and $\ftt$ in $\overline{\Omega}_{2}$) such that there exists a finite number of points $x_{1,1}^{j}\in\widetilde{\Omega}_{1}$ for $j=1,\dots,N_{1,1}$ and $x_{1,2}^{j}\in\widetilde{\Omega}_{1}$ for $j=1,\dots,N_{1,2}$ (resp. $x_{2,1}^{j}\in\Omega_{2}$ for $j=1,\dots,N_{2,1}$ and $x_{2,2}^{j}\in\Omega_{2}$ for $j=1,\dots,N_{2,2}$) and $\epsilon>0$ such that $\overline{B(x_{1,1}^{j},2\epsilon)}\subset\widetilde{\Omega}_{1}$, $\overline{B(x_{1,2}^{j},2\epsilon)}\subset\widetilde{\Omega}_{1}$, $B(x_{1,1}^{j_{1}},2\epsilon)\cap B(x_{1,2}^{j_{2}},2\epsilon)=\emptyset$ and in $B(x_{1,k}^{j},2\epsilon)$ we have $\vfi_{1,k+1}>\vfi_{1,k}$ (resp. $\overline{B}(x_{2,1}^{j},2\epsilon)\subset\Omega_{2}$, $\overline{B}(x_{2,2}^{j},2\epsilon)\subset\Omega_{2}$, $B(x_{2,1}^{j_{1}},2\epsilon)\cap B(x_{1,2}^{j_{2}},2\epsilon)=\emptyset$ and in $B(x_{2,k}^{j},2\epsilon)$ we have $\vfi_{2,k+1}>\vfi_{2,k}$), where $k+1$ is equal to $2$ if $k=1$ and equal to $1$ if $k=2$. Furthermore, by setting $\displaystyle U_{1,k}=\widetilde{\Omega}_{1}\bigcap\left(\bigcup_{j=1}^{N_{1,k}}B(x_{1,k}^{j},2\epsilon)\right)^{c}$ and $\displaystyle U_{2,k}=\Omega_{2}\bigcap\left(\bigcup_{j=1}^{N_{2,k}}B(x_{2,k}^{j},2\epsilon)\right)^{c}$ for $k=1,2$, $\gamma_{1}=\p B_{r}$, $\gamma_{2}=\Gamma_{2}$ and $\gamma=S$, the phases verifying that $|\nabla\foo|>0$ in $\overline{U}_{1,1}$, $|\nabla\phot|>0$ in $\overline{U}_{1,2}$,$|\nabla\phto|>0$ in $\overline{U}_{2,1}$, $|\nabla\ftt|>0$ in $\overline{U}_{2,2}$, $\pn\vfi_{1,k\,|\gamma_{1}}\neq0$, $\pn\vfi_{2,k\,|\gamma_{2}}<0$, $\pn\vfi_{1,k\,|\gamma}< 0$, $\pn\vfi_{2,k\,|\gamma}< 0$, and $\foo,\,\phot,\,\phto$ and $\ftt$ verifying the H{\"o}rmander's condition~\eqref{c5} respectively in $U_{1,1},\, U_{1,2},\,U_{2,1}$ and $U_{2,2}$. We can let also  (see~\cite{B}) $\vfi_{1,1\,|S}=\vfi_{2,1\,|S}$ where by construction we obtain $\vfi_{1,2|S}=\vfi_{2,2|S}$  and by argument of density we can suppose that $(\pn\vfi_{1,1\,|S})^{2}-(\pn\vfi_{2,1\,|S})^{2}>1$ and $(\pn\vfi_{1,2\,|S})^{2}-(\pn\vfi_{2,2\,|S})^{2}>1$. And that concludes our construction of the weight functions.
\section{Resolvent Estimate}\label{f1}
This section is devoted to prove a resolvent estimate, precisely we will show that for some constant $C>0$ we have
\begin{equation}\label{f2}
\|(\mathcal{A}-i\mu\,\id)^{-1}\|_{\mathcal{L}(\mathcal{H})}\leq C\e^{C|\mu|},
\end{equation}
for every $\mu\in\R$ large enough in absolute value, Which by Burq's result follow the kind of energy decay rate given in Theorem~\ref{b2}.

We suppose that the resolvent estimate~\eqref{f2} is false. Then there exist two sequences $K_{m}>0$ and $\mu_{m}\in\R$ and two families $(u_{1,m},u_{2,m},v_{1,m},v_{2,m})\in\mathcal{D}(\mathcal{A})$ and $(f_{1,m},f_{2,m},g_{1,m},g_{2,m})\in\mathcal{H}$, $m=1,2,\ldots$  such that
\begin{equation}\label{f21}
|\mu_{m}|\,\longrightarrow\,+\infty,\qquad K_{m}\,\longrightarrow\,+\infty,\qquad\|(u_{1,m},u_{2,m},v_{1,m},v_{2,m})\|_{\mathcal{H}}=1,
\end{equation}
and
\begin{equation}\label{f22}
\e^{K_{m}|\mu_{m}|}(\mathcal{A}-i\mu_{m})\left(\begin{array}{c}
u_{1,m}
\\
u_{2,m}
\\
v_{1,m}
\\
v_{2,m}
\end{array}\right)=\left(\begin{array}{c}
f_{1,m}
\\
f_{2,m}
\\
g_{1,m}
\\
g_{2,m}
\end{array}\right)\,\longrightarrow\,0\text{ in }\mathcal{H}.
\end{equation}
This implies that
\begin{equation}\label{f23}
\e^{K_{m}|\mu_{m}|}\left(\begin{array}{c}
v_{1,m}-i\mu_{m}u_{1,m}
\\
v_{2,m}-i\mu_{m}u_{2,m}
\\
-\Delta(\Delta u_{1,m}+a\Delta v_{1,m})-i\mu_{m}v_{1,m}
\\
\Delta u_{2,m}-i\mu_{m}v_{2,m}
\end{array}\right)=\left(\begin{array}{c}
f_{1,m}
\\
f_{2,m}
\\
g_{1,m}
\\
g_{2,m}
\end{array}\right)\,\longrightarrow\,0\text{ in }\mathcal{H}.
\end{equation}
From~\eqref{f21} and~\eqref{f22}, we get
\begin{equation}\label{f25}
\re\left\langle\left(\begin{array}{c}
f_{1,m}
\\
f_{2,m}
\\
g_{1,m}
\\
g_{2,m}
\end{array}\right),\left(\begin{array}{c}
u_{1,m}
\\
u_{2,m}
\\
v_{1,m}
\\
v_{2,m}
\end{array}\right)\right\rangle_{\mathcal{H}}=-\e^{K_{m}|\mu_{m}|}\int_{\Omega_{1}}a|\Delta v_{1,m}|^{2}\,\ud x\,\longrightarrow\,0.
\end{equation}
Then from the first equation of~\eqref{f23} and~\eqref{f25}, we obtain
\begin{equation}\label{f26}
|\mu_{m}|^{2}\e^{K_{m}|\mu_{m}|}\int_{\Omega_{1}}a|\Delta u_{1,m}|^{2}\,\ud x\,\longrightarrow\,0.
\end{equation}
Since, from~\eqref{f23} we have
$$
\Delta v_{1,m}=i\mu_{m}\e^{-K_{m}|\mu_{m}|}\Delta u_{1,m}+\e^{-K_{m}|\mu_{m}|}\Delta f_{1,m}
$$
then by elliptic estimates it follows that
$$
\|v_{1,m}\|_{H^{2}(\omega)}\leq C\left(|\mu_{m}|.\|\Delta u_{1,m}\|_{L^{2}(\Omega_{1})}\e^{-K_{m}|\mu_{m}|}+\|\Delta f_{1,m}\|_{L^{2}(\Omega_{1})}\e^{-K_{m}|\mu_{m}|}+\|v_{1,m}\|_{L^{2}(\Omega_{1})}\right)
$$
which mean that $\|v_{1,m}\|_{H^{2}(\omega)}$ is bounded. We multiply the third equation of~\eqref{f23} by $\mu_{m}^{-1}\psi.\overline{v}_{1,m}$ where $\psi\in\Ci(\Omega_{1})$ and $\supp(\psi)\subset\omega$ then from~\eqref{f25} and~\eqref{f26} we obtain
\begin{equation*}
\e^{\frac{K_{m}}{4}|\mu_{m}|}\int_{\omega}|v_{1,m}|^{2}\psi\,\ud x\,\longrightarrow\,0.
\end{equation*}
In particular we have
\begin{equation*}
\e^{\frac{K_{m}}{4}|\mu_{m}|}\int_{B_{4r}}|v_{1,m}|^{2}\,\ud x\,\longrightarrow\,0.
\end{equation*}
So that, from the first equation of~\eqref{f23} we show that
\begin{equation}\label{f29}
\e^{\frac{K_{m}}{4}|\mu_{m}|}\int_{B_{4r}}|u_{1,m}|^{2}\,\ud x\,\longrightarrow\,0.
\end{equation}

We consider now the following transmission problem
\begin{equation}\label{f3}
\left\{\begin{array}{ll}
v_{1}-i\mu u_{1}=f_{1}&\text{in }\Omega_{1}
\\
v_{2}-i\mu u_{2}=f_{2}&\text{in }\Omega_{2}
\\
-\Delta(\Delta u_{1}+a\Delta v_{1})-i\mu v_{1}=g_{1}&\text{in }\Omega_{1}
\\
\Delta u_{2}-i\mu v_{2}=g_{2}&\text{in }\Omega_{2}
\\
u_{1}=u_{2},\quad\pn u_{1}=0&\text{on }S
\\
\pn\Delta u_{1}+\pn u_{2}=0&\text{on }S
\\
u_{2}=0&\text{on }\Gamma.
\end{array}\right.
\end{equation}
By setting
\begin{equation}\label{f4}
\left\{\begin{array}{ll}
\Phi_{1}=g_{1}+i\mu f_{1}
\\
\Phi_{2}=g_{2}+i\mu f_{2},
\end{array}\right.
\end{equation}
then~\eqref{f3} can be recast as follows
\begin{equation}\label{f5}
\left\{\begin{array}{ll}
v_{1}=f_{1}+i\mu u_{1}&\text{in }\Omega_{1}
\\
v_{2}=f_{2}+i\mu u_{2}&\text{in }\Omega_{2}
\\
-\Delta(\Delta u_{1}+a\Delta v_{1})+\mu^{2} u_{1}=\Phi_{1}&\text{in }\Omega_{1}
\\
-\Delta u_{2}-\mu^{2} u_{2}=-\Phi_{2}&\text{in }\Omega_{2}
\\
u_{1}=u_{2},\quad\pn u_{1}=0&\text{on }S
\\
\pn\Delta u_{1}+\pn u_{2}=0&\text{on }S
\\
u_{2}=0&\text{on }\Gamma.
\end{array}\right.
\end{equation}
We denote by
\begin{equation}\label{f6}
\left\{\begin{array}{l}
z_{1}'=\Delta u_{1}+a\Delta v_{1}-|\mu|u_{1}
\\
z_{1}''=\Delta u_{1}+a\Delta v_{1}
\\
z_{2}=-u_{2},
\end{array}\right.
\end{equation}
it follows from~\eqref{f5} that $z_{1}',\,z_{1}''$ and $z_{2}$ are solution of the following transmission problem
\begin{equation}\label{f7}
\left\{\begin{array}{ll}
-\Delta z_{1}'-|\mu|z_{1}'=\Phi_{1}-|\mu|a\Delta v_{1}&\text{in }\Omega_{1}
\\
-\Delta z_{1}''+|\mu|z_{1}''=\Phi_{1}+|\mu|z_{1}'&\text{in }\Omega_{1}
\\
-\Delta z_{2}-|\mu|^{2} z_{2}=\Phi_{2}&\text{in }\Omega_{2}
\\
\frac{1}{|\mu|}z_{1}'=z_{2}+\frac{1}{|\mu|}z_{1}''&\text{on }S
\\
\pn z_{1}'=\pn z_{2}&\text{on }S
\\
\frac{1}{|\mu|}z_{1}''=z_{2}+\frac{1}{|\mu|}\theta&\text{on }S
\\
\pn z_{1}''=\pn z_{2}&\text{on }S
\\
z_{2}=0&\text{on }\Gamma,
\end{array}\right.
\end{equation}
where $\theta=-z_{1}'+2z_{1}''$.

We set $B_{5r}$ a ball of radius $5r>0$ such that $a(x)>0$ in $B_{5r}\subset\omega$. We set as the previous section $\widetilde{\Omega}_{1}=\Omega_{1}\backslash \overline{B}_{r}$. The most important ingredient of the proof of the resolvent estimate~\eqref{f2} is the following lemma which is essentially the result of the Carleman estimate.
\begin{lem}
There exist a constant $C>0$ such that for any $(u_{1},u_{2},v_{1},v_{2})\in\mathcal{D}(\mathcal{A})$ solution of~\eqref{f3} the following estimate holds
\begin{equation}\label{f8}
\begin{split}
\|u_{1}\|_{L^{2}(\Omega_{1})}^{2}+\|\Delta u_{1}\|_{L^{2}(\Omega_{1})}^{2}+\|u_{2}\|_{L^{2}(\Omega_{2})}^{2}+\|\nabla u_{2}\|_{L^{2}(\Omega_{2})}^{2}\leq C\e^{C/h}\bigg(\|\Delta f_{1}\|_{L^{2}(\Omega_{1})}^{2}
\\
+\|f_{2}\|_{L^{2}(\Omega_{2})}^{2}+\|g_{1}\|_{L^{2}(\Omega_{1})}^{2}+\|g_{2}\|_{L^{2}(\Omega_{2})}^{2}+\int_{\Omega_{1}}a|\Delta u_{1}|^{2}\,\ud x+\int_{B_{5r}}|u_{1}|^{2}\,\ud x\bigg).
\end{split}
\end{equation}
for all $|\mu|>0$ large enough.
\end{lem}
\begin{pr}
We introduce the cut-off function $\chi\in\Ci(\Omega_{1})$ by setting
$$
\chi(x)=\left\{\begin{array}{ll}
1&\text{in }B_{3r}^{c}
\\
0&\text{in }B_{2r}.
\end{array}\right.
$$
Next, we denote by $\tilde{z}_{1}'=\chi z_{1}'$ and $\tilde{z}_{1}''=\chi z_{1}''$. Then by~\eqref{f7}, one sees that
\begin{equation}\label{f9}
\left\{\begin{array}{l}
-\Delta\tilde{z}_{1}'-|\mu|\tilde{z}_{1}'=\widetilde{\Phi}_{1}'=\chi\Phi_{1}-|\mu|\chi a\Delta v_{1}-[\Delta,\chi]z_{1}'
\\
-\Delta\tilde{z}_{1}''+|\mu|\tilde{z}_{1}''=\widetilde{\Phi}_{1}''=\chi\Phi_{1}+|\mu|\tilde{z}_{1}'-[\Delta,\chi]z_{1}''.
\end{array}\right.
\end{equation}
Keeping the same notations as the end of the previous section, and focus now to the system~\eqref{f7}. Taking $\foo$, $\phot$, $\phto$, and $\ftt$ the four weight functions that satisfy the conclusion of the end of section~\ref{b1}. We set $\chi_{1,1}$, $\chi_{1,2}$, $\chi_{2,1}$ and $\chi_{2,2}$ four cut-off functions that equal to one  respectively  in $\displaystyle\left(\bigcup_{j=1}^{N_{1,1}}B(x_{1,1}^{j},2\epsilon)\right)^{c}$, $\displaystyle\left(\bigcup_{j=1}^{N_{1,2}}B(x_{1,2}^{j},2\epsilon)\right)^{c}$, $\displaystyle\left(\bigcup_{j=1}^{N_{2,1}}B(x_{2,1}^{j},2\epsilon)\right)^{c}$ and $\displaystyle\left(\bigcup_{j=1}^{N_{2,2}}B(x_{2,2}^{j},2\epsilon)\right)^{c}$ and supported in $\displaystyle\left(\bigcup_{j=1}^{N_{1,1}}B(x_{1,1}^{j},\epsilon)\right)^{c}$, $\displaystyle\left(\bigcup_{j=1}^{N_{1,2}}B(x_{1,2}^{j},\epsilon)\right)^{c}$, $\displaystyle\left(\bigcup_{j=1}^{N_{2,1}}B(x_{2,1}^{j},\epsilon)\right)^{c}$ and $\displaystyle\left(\bigcup_{j=1}^{N_{2,2}}B(x_{2,2}^{j},\epsilon)\right)^{c}$ respectively (in order to eliminate the critical points of the phases functions $\foo$, $\phot$, $\phto$, $\ftt$, $\phto$ and $\ftt$). We set $y_{1,1}'=\chi_{1,1}\tilde{z}_{1}'$, $y_{1,2}'=\chi_{1,2}\tilde{z}_{1}'$, $y_{1,1}''=\chi_{1,1}\tilde{z}_{1}''$, $y_{1,2}''=\chi_{1,2}\tilde{z}_{1}''$, $y_{2,1}=\chi_{2,1}z_{2}$ and $y_{2,2}=\chi_{2,2}z_{2}$. Then from~\eqref{f7} and~\eqref{f9} and by noting $\displaystyle h=\frac{1}{|\mu|}$  for $k=1,2$ we obtain
\begin{equation}\label{f10}
\left\{\begin{array}{ll}
\displaystyle-\Delta y_{1,k}'-y_{1,k}'/h=\Psi_{1,k}'&\text{in }U_{1,k}
\\
\displaystyle-\Delta y_{1,k}''+y_{1,k}''/h=\Psi_{1,k}''&\text{in }U_{1,k}
\\
\displaystyle-\Delta y_{2,k}-y_{2,k}/h^{2}=\Psi_{2,k}&\text{in }U_{2,k}
\\
hy_{1,k}'=y_{2,k}+hy_{1,k}''&\text{on }S
\\
\pn y_{1,k}'=\pn y_{2,k}&\text{on }S
\\
hy_{1,k}''=y_{2,k}+h\theta&\text{on }S
\\
\pn y_{1,k}''=\pn y_{2,k}&\text{on }S
\\
y_{2,k}=0&\text{on }\Gamma,
\end{array}\right.
\end{equation}
where
\begin{equation}\label{f11}
\left\{\begin{array}{l}
\Psi_{1,k}'=-[\Delta,\chi_{1,k}]z_{1}'+\chi_{1,k}\widetilde{\Phi}_{1}'
\\
\Psi_{1,k}''=-[\Delta,\chi_{1,k}]z_{1}''+\chi_{1,k}\widetilde{\Phi}_{1}''
\\
\Psi_{2,k}=-[\Delta,\chi_{2,k}]z_{2}+\chi_{2,k}\Phi_{2}.
\end{array}\right.
\end{equation}
Applying Carleman estimate of Theorem~\ref{c4} to the systems~\eqref{f10} then for $k=1,2$ we obtain
\begin{equation}\label{f12}
\begin{split}
h\|\e^{\vfi_{1,k}/h}y_{1,k}'\|_{L^{2}(U_{1,k})}^{2}+h\|\e^{\vfi_{1,k}/h}y_{1,k}''\|_{L^{2}(U_{1,k})}^{2}+h\|\e^{\vfi_{2,k}/h}y_{2,k}\|_{L^{2}(U_{2,k})}^{2}+
\\
h^{3}\|\e^{\vfi_{1,k}/h}\nabla y_{1,k}'\|_{L^{2}(U_{1,k})}^{2}+h^{3}\|\e^{\vfi_{1,k}/h}\nabla y_{1,k}''\|_{L^{2}(U_{1,k})}^{2}+h^{3}\|\e^{\vfi_{2,k}/h}\nabla y_{2,k}\|_{L^{2}(U_{2,k})}^{2}+
\\
h|\e^{\vfi_{1,k}/h}y_{1,k}'|_{L^{2}(S)}^{2}+h|\e^{\vfi_{1,k}/h}\nabla y_{1,k}'|_{L^{2}(S)}^{2}+h|\e^{\vfi_{1,k}/h}y_{1,k}''|_{L^{2}(S)}^{2}+
\\
h|\e^{\vfi_{1,k}/h}\nabla y_{1,k}''|_{L^{2}(S)}^{2}\leq C(h^{4}\|\e^{\vfi_{1,k}/h}\Psi_{1,k}'\|_{L^{2}(U_{1,k})}^{2}+h^{4}\|\e^{\vfi_{1,k}/h}\Psi_{1,k}''\|_{L^{2}(U_{1,k})}^{2}
\\
+h^{4}\|\e^{\vfi_{2,k}/h}\Psi_{2,k}\|_{L^{2}(U_{2,k})}^{2}+h^{3}|\e^{\vfi_{1,k}/h}\theta|_{L^{2}(S)}^{2}+h^{5}|\e^{\vfi_{1,k}/h}\nabla \theta|_{L^{2}(S)}^{2}).
\end{split}
\end{equation}
The two last terms of the right hand side of~\eqref{f12} can be absorbed to the left hand side for $h>0$ small enough and since $\theta=-y_{1,k}'+2y_{1,k}''$, therefore by~\eqref{f11} we arrive at
\begin{equation}\label{f13}
\begin{split}
h\|\e^{\vfi_{1,k}/h}y_{1,k}'\|_{L^{2}(U_{1,k})}^{2}+h\|\e^{\vfi_{1,k}/h}y_{1,k}''\|_{L^{2}(U_{1,k})}^{2}+h\|\e^{\vfi_{2,k}/h}y_{2,k}\|_{L^{2}(U_{2,k})}^{2}
\\
+h^{3}\|\e^{\vfi_{1,k}/h}\nabla y_{1,k}'\|_{L^{2}(U_{1,k})}^{2}+h^{3}\|\e^{\vfi_{1,k}/h}\nabla y_{1,k}''\|_{L^{2}(U_{1,k})}^{2}+h^{3}\|\e^{\vfi_{2,k}/h}\nabla y_{2,k}\|_{L^{2}(U_{2,k})}^{2}
\\
\leq Ch^{4}(\|\e^{\vfi_{1,k}/h}\widetilde{\Phi}_{1}'\|_{L^{2}(U_{1,k})}^{2}+\|\e^{\vfi_{1,k}/h}\widetilde{\Phi}_{1}''\|_{L^{2}(U_{1,k})}^{2}+\|\e^{\vfi_{2,k}/h}\Phi_{2}\|_{L^{2}(U_{2,k})}^{2}
\\
+\|\e^{\vfi_{1,k}}[\Delta,\chi_{1,k}]\tilde{z}_{1}'\|_{L^{2}(U_{1,k})}^{2}+\|\e^{\vfi_{1,k}}[\Delta,\chi_{1,k}]\tilde{z}_{1}''\|_{L^{2}(U_{1,k})}^{2}+\|\e^{\vfi_{2,k}}[\Delta,\chi_{2,k}]z_{2}\|_{L^{2}(U_{2,k})}^{2}).
\end{split}
\end{equation}
We addition the two last estimates  for $k=1,2$ and using the properties of phases $\vfi_{1,k}<\vfi_{1,k+1}$ in $\displaystyle\left(\bigcup_{j=1}^{N_{1,k}}B(x_{1,k}^{j},2\epsilon)\right)$ and $\vfi_{2,k}<\vfi_{2,k+1}$ in $\displaystyle\left(\bigcup_{j=1}^{N_{2,k}}B(x_{2,k}^{j},2\epsilon)\right)$ then we can absorb the terms $[\Delta,\chi_{1,k}]\tilde{z}_{1}'$, $[\Delta,\chi_{1,k}]\tilde{z}_{1}''$ and $[\Delta,\chi_{2,k}]z_{2}$ at the right hand side of~\eqref{f13} into the left hand side for $h>0$ small. Namely, we find
\begin{equation}\label{f14}
\begin{split}
h\int_{\widetilde{\Omega}_{1}}(\e^{2\foo/h}+\e^{2\phot/h})|\tilde{z}_{1}'|^{2}\,\ud x+h\int_{\widetilde{\Omega}_{1}}(\e^{2\foo/h}+\e^{2\phot/h})|\tilde{z}_{1}''|^{2}\,\ud x
\\
+h\int_{\Omega_{2}}(\e^{2\phto/h}+\e^{2\ftt/h})|z_{2}|^{2}\,\ud x+h\int_{\Omega_{2}}(\e^{2\phto/h}+\e^{2\ftt/h})|\nabla z_{2}|^{2}\,\ud x 
\\
\leq C\bigg(h^{4}\int_{\Omega_{1}}(\e^{2\foo/h}+\e^{2\phot/h})|\Phi_{1}|^{2}\,\ud x+h^{4}\int_{\Omega_{2}}(\e^{2\phto/h}+\e^{2\ftt/h})|\Phi_{2}|^{2}\,\ud x
\\
h^{4}\int_{\widetilde{\Omega}_{1}}(\e^{2\foo/h}+\e^{2\phot/h})|[\Delta,\chi]z_{1}'|^{2}\,\ud x+h^{4}\int_{\widetilde{\Omega}_{1}}(\e^{2\foo/h}+\e^{2\phot/h})|[\Delta,\chi]z_{1}''|^{2}\,\ud x
\\
+h^{2}\int_{\widetilde{\Omega}_{1}}(\e^{2\foo/h}+\e^{2\phot/h})|a\Delta v_{1}|^{2}\,\ud x+h^{2}\int_{\widetilde{\Omega}_{1}}(\e^{2\foo/h}+\e^{2\phot/h})|z_{1}'|^{2}\,\ud x\bigg).
\end{split}
\end{equation}
For $h>0$ small we can absorb the last term of the right hand side of~\eqref{f14} into the left hand side. Besides, by remarking that $\Omega_{1}=\widetilde{\Omega}_{1}\cup B_{3r}$ and by taking the maximum of $\foo$, $\phot$, $\phto$ and $\ftt$ into the right hand side of~\eqref{f14} and their minimum into the left hand side then it follows from the definitions of $\Phi_{1}$ and $\Phi_{2}$ in~\eqref{f4} that
\begin{equation}\label{f16}
\begin{split}
\int_{\Omega_{1}}|z_{1}'|^{2}\,\ud x+\int_{\Omega_{1}}|z_{1}''|^{2}\,\ud x+\int_{\Omega_{2}}|z_{2}|^{2}\,\ud x+\int_{\Omega_{2}}|\nabla z_{2}|^{2}\,\ud x\leq C\e^{C/h}\bigg(\int_{\Omega_{1}}|f_{1}|^{2}\,\ud x
\\
+\int_{\Omega_{2}}|f_{2}|^{2}\,\ud x+\int_{\Omega_{1}}|g_{1}|^{2}\,\ud x+\int_{\Omega_{2}}|g_{2}|^{2}\,\ud x+\int_{B_{3r}}|z_{1}'|^{2}\,\ud x+\int_{B_{3r}}|z_{1}''|^{2}\,\ud x
\\
+\int_{\Omega_{1}}a|\Delta v_{1}|^{2}\,\ud x+\int_{\Omega_{1}}|[\Delta,\chi]z_{1}'|^{2}\,\ud x+\int_{\Omega_{1}}|[\Delta,\chi]z_{1}''|^{2}\,\ud x\bigg).
\end{split}
\end{equation}
Let $\widetilde{\chi}$ be a cut-off function equal to $1$ in a neighborhood of $B_{4r}$ and supported in $B_{5r}$ then by the second equation of~\eqref{f5} and of~\eqref{f6} we have
$$
(-1+\Delta)(\widetilde{\chi}z_{1}'')=[\Delta,\widetilde{\chi}]z_{1}''-\widetilde{\chi}z_{1}''-|\mu|^{2}\widetilde{\chi}u_{1}-\widetilde{\chi}g_{1}-i\mu\widetilde{\chi} f_{1}.
$$
Hence by elliptic estimates (see~\cite{WRL}) we get
\begin{eqnarray}\label{f17}
\|z_{1}''\|_{H^{1}(B_{4r})}^{2}&\leq& C(\|(-1+\Delta)(\widetilde{\chi}z_{1}'')\|_{H^{-1}(B_{5r})}^{2}+\|z_{1}''\|_{L^{2}(B_{5r})}^{2})\nonumber
\\
&\leq& C(\|g_{1}\|_{L^{2}(\Omega_{1})}^{2}+|\mu|^{2}\|f_{1}\|_{L^{2}(\Omega_{1})}^{2}+|\mu|^{4}\|u_{1}\|_{L^{2}(B_{5r})}^{2}+\|z_{1}''\|_{L^{2}(B_{5r})}^{2}).
\end{eqnarray}
Since $\supp([\Delta,\chi])\subset B_{3r}$ we deduce from~\eqref{f6} and~\eqref{f17} that
\begin{equation}\label{f18}
\begin{split}
\int_{B_{3r}}|z_{1}''|^{2}\,\ud x+\int_{\Omega_{1}}|[\Delta,\chi]z_{1}''|^{2}\,\ud x\leq C\|z_{1}''\|_{H^{1}(B_{4r})}^{2}\leq C\Big(\|g_{1}\|_{L^{2}(\Omega_{1})}^{2}+|\mu|^{2}\|f_{1}\|_{L^{2}(\Omega_{1})}^{2}
\\
+\|\Delta u_{1}\|_{L^{2}(B_{5r})}^{2}+\|a\Delta v_{1}\|_{L^{2}(B_{5r})}^{2}+|\mu|^{4}\|u_{1}\|_{L^{2}(B_{5r})}^{2}\Big).
\end{split}
\end{equation}
Similarly, we prove also that
\begin{equation}\label{f15}
\begin{split}
\int_{B_{3r}}|z_{1}'|^{2}\,\ud x+\int_{\Omega_{1}}|[\Delta,\chi]z_{1}'|^{2}\,\ud x\leq C\Big(\|g_{1}\|_{L^{2}(\Omega_{1})}^{2}+|\mu|^{2}\|f_{1}\|_{L^{2}(\Omega_{1})}^{2}
\\
+\|\Delta u_{1}\|_{L^{2}(B_{5r})}^{2}+\|a\Delta v_{1}\|_{L^{2}(B_{5r})}^{2}+|\mu|^{2}\|u_{1}\|_{L^{2}(B_{5r})}^{2}\Big).
\end{split}
\end{equation}
We combine~\eqref{f16} with~\eqref{f18} and~\eqref{f15} and we recall the expression of $z_{1}''$ and $z_{2}$ in~\eqref{f6} then we find
\begin{equation}\label{f19}
\begin{split}
\|\Delta u_{1}\|_{L^{2}(\Omega_{1})}^{2}+\|\nabla u_{2}\|_{L^{2}(\Omega_{2})}^{2}+\|u_{2}\|_{L^{2}(\Omega_{2})}^{2}\leq C\e^{C/h}\bigg(\|f_{1}\|_{L^{2}(\Omega_{1})}^{2}+\|f_{2}\|_{L^{2}(\Omega_{2})}^{2}+
\\
\|g_{1}\|_{L^{2}(\Omega_{1})}^{2}+\|g_{2}\|_{L^{2}(\Omega_{2})}^{2}+\int_{B_{5r}}|\Delta u_{1}|^{2}\,\ud x+\int_{\Omega_{1}}a|\Delta v_{1}|^{2}\,\ud x+\int_{B_{5r}}|u_{1}|^{2}\,\ud x\bigg).
\end{split}
\end{equation}
We substitute the expression of $v_{1}$ and $v_{2}$ in~\eqref{f5} into~\eqref{f19} then we obtain
\begin{equation}\label{f20}
\begin{split}
\|\Delta u_{1}\|_{L^{2}(\Omega_{1})}^{2}+\|\nabla u_{2}\|_{L^{2}(\Omega_{2})}^{2}+\|u_{2}\|_{L^{2}(\Omega_{2})}^{2}\leq C\e^{C/h}\bigg(\|f_{1}\|_{L^{2}(\Omega_{1})}^{2}+\|\Delta f_{1}\|_{L^{2}(\Omega_{1})}^{2}
\\
+\|f_{2}\|_{L^{2}(\Omega_{2})}^{2}+\|g_{1}\|_{L^{2}(\Omega_{1})}^{2}+\|g_{2}\|_{L^{2}(\Omega_{2})}^{2}+\int_{\Omega_{1}}a|\Delta u_{1}|^{2}\,\ud x+\int_{B_{5r}}|u_{1}|^{2}\,\ud x\bigg).
\end{split}
\end{equation}
The estimate~\eqref{f8} holds now from~\eqref{f20}, Poincar\'e inequality and Lemma~\ref{b5}.
\end{pr}

Applying inequality~\eqref{f8} to the system~\eqref{f23} it follows that
\begin{equation*}
\begin{split}
\|\Delta u_{1,m}\|_{L^{2}(\Omega_{1})}^{2}+\|\nabla u_{2,m}\|_{L^{2}(\Omega_{2})}^{2}+\|u_{1,m}\|_{L^{2}(\Omega_{1})}^{2}+\|u_{2,m}\|_{L^{2}(\Omega_{2})}^{2}\leq
\\
C\e^{C|\mu_{m}|}\bigg(\e^{-2K_{m}|\mu_{m}|}\Big(\|\Delta f_{1,m}\|_{L^{2}(\Omega_{1})}^{2}+\|\nabla f_{2,m}\|_{L^{2}(\Omega_{2})}^{2}+\|g_{1,m}\|_{L^{2}(\Omega_{1})}^{2}+\|g_{2,m}\|_{L^{2}(\Omega_{2})}^{2}\Big)
\\
+\left(\int_{\Omega_{1}}a|\Delta u_{1,m}|^{2}\,\ud x+\int_{B_{4r}}|u_{1,m}|^{2}\,\ud x\right)\bigg).
\end{split}
\end{equation*}
We use the expression of $u_{1,m}$ and $u_{2,m}$ in~\eqref{f23} we follows that
\begin{equation}\label{f31}
\begin{split}
\|\Delta u_{1,m}\|_{L^{2}(\Omega_{1})}^{2}+\|\nabla u_{2,m}\|_{L^{2}(\Omega_{2})}^{2}+\|v_{1,m}\|_{L^{2}(\Omega_{1})}^{2}+\|v_{2,m}\|_{L^{2}(\Omega_{2})}^{2}\leq
\\
C\e^{C|\mu_{m}|}\bigg(\e^{-2K_{m}|\mu_{m}|}\Big(\|\Delta f_{1,m}\|_{L^{2}(\Omega_{1})}^{2}+\|\nabla f_{2,m}\|_{L^{2}(\Omega_{2})}^{2}+\|g_{1,m}\|_{L^{2}(\Omega_{1})}^{2}+\|g_{2,m}\|_{L^{2}(\Omega_{2})}^{2}\Big)
\\
+\left(\int_{\Omega_{1}}a|\Delta u_{1,m}|^{2}\,\ud x+\int_{B_{4r}}|u_{1,m}|^{2}\,\ud x\right)\bigg).
\end{split}
\end{equation}
Finally~\eqref{f21},~\eqref{f22} and~\eqref{f25} and~\eqref{f29} shows that the right hand side of~\eqref{f31} go to zero as $m\,\longrightarrow\,+\infty$, hence we obtain a contradiction with~\eqref{f21}, therefore the resolvent estimate~\eqref{f2} is proved now.

Now, follows to~\cite[Lemma 4.1]{CLL} it just remains to show that $\mathcal{A}$ has no purely imaginary eigenvalue. Further, $0\in\rho(\mathcal{A})$, where $\rho(\mathcal{A})$ stands for the resolvent set of $\mathcal{A}$. Let $\mu\neq 0$ be a real number. Suppose that for some $U=(u_{1},u_{2},v_{1},v_{2})\in\mathcal{D}(\mathcal{A})$, one has
\begin{equation}\label{f32}
\mathcal{A}U=i\mu U.
\end{equation}
We shall show that $U=0$. Taking the inner product with $U$ on both side of~\eqref{f32} and taking the real part we immediately find that $v_{1}=0$ in $\supp(a)$. Now~\eqref{f32} can be recast as
\begin{equation}\label{f33}
\left\{\begin{array}{ll}
v_{1}=i\mu_{1}u_{1}&\text{in }\Omega_{1},
\\
v_{2}=i\mu_{2}u_{2}&\text{in }\Omega_{2},
\\
-\Delta(\Delta u_{1}+a\Delta v_{1})-i\mu v_{1}=0&\text{in }\Omega_{1},
\\
\Delta u_{2}-i\mu v_{2}=0&\text{in }\Omega_{2}.
\end{array}\right.
\end{equation}
Since $v_{1}=0$ in $\supp(a)$ and $\mu\neq 0$ the top line of~\eqref{f33} yields $u_{1}=0$ in $\supp(a)$. The third line of~\eqref{f33} combined with the first one could be written as
$$
\Delta z+\mu z=0\quad \text{in }\Omega_{1}\qquad\text{and}\qquad z=0\quad \text{in }\supp(a),
$$
where we denoted by $z=\Delta u_{1}-\mu u_{1}$. Since $\omega\subset\supp(a)$ then by Calder\'on's theorem~\cite[Theorem 4.2]{RR2} for elliptic operators we find that $z=0$, this mean that $\Delta u_{1}-\mu u_{1}=0$ which imply for the same argument as previously that $u_{1}=0$ in $\Omega_{1}$. Reporting that in the first line of~\eqref{f33}, we derive $v_{1}=0$ in $\Omega_{1}$. The second and fourth line of~\eqref{f33} lead to
$$
\Delta u_{2}+\mu^{2}u_{2}=0\quad\text{in }\Omega_{2}
$$
with the boundary conditions
$$
u_{2}=0\quad\text{on }\p\Omega_{2}\quad\text{and}\quad \pn u_{2}=0\quad\text{on } S.
$$
By standard theory in linear elliptic equations $u_{2}=0$ in $\Omega_{2}$. Using the second line of~\eqref{f33}, we get $v_{2}=0$ in $\Omega_{2}$; hence $U=0$. Therefore, $\mathcal{A}$ has no purely eigenvalue.
\subsubsection*{Acknowledgments}
The author thanks the referees for many valuable remarks which helped us to improve the paper significantly. 
\nocite{*}
\bibliographystyle{alpha}
\bibliography{bibPWKV}
\addcontentsline{toc}{section}{References}
\end{document}